# RADIAL SYMMETRY OF $p$-HARMONIC MINIMIZERS

ALEKSIS KOSKI AND JANI ONNINEN

ABSTRACT. "It is still not known if the radial cavitating minimizers obtained by Ball [J.M. Ball, Discontinuous equilibrium solutions and cavitation in nonlinear elasticity, Phil. Trans. R. Soc. Lond. A 306 (1982) 557–611] (and subsequently by many others) are global minimizers of any physically reasonable nonlinearly elastic energy". The quotation is from [37] and seems to be still accurate. The model case of the $p$-harmonic energy is considered here. We prove that the planar radial minimizers are indeed the global minimizers provided we prescribe the admissible deformations on the boundary. In the traction free setting, however, even the identity map need not be a global minimizer.

## 1. INTRODUCTION

The general law of hyperelasticity tells us that there exists an energy integral functional with a given *stored-energy* function characterizing the elastic properties of a material. The mathematical models of nonlinear elasticity have been pioneered by Antman [2], Ball [5], and Ciarlet [13]. One of important problems in nonlinear elasticity is whether or not the radially symmetric minimizers are indeed global minimizers of the given physically reasonable energy. This leads us to study energy-minimal homeomorphisms $h \colon \mathbb{A} \xrightarrow{\text{onto}} \mathbb{A}^*$ between annuli

$$\mathbb{A} = A(r,R) = \{z \in \mathbb{R}^n : r < |z| < R\} \quad \text{and} \quad \mathbb{A}^* = \{z \in \mathbb{R}^n : r_* < |z| < R_*\}.$$

Hereafter $0 \leqslant r < R$ and $0 \leqslant r_* < R_*$ are called the inner and outer radii of $\mathbb{A}$ and $\mathbb{A}^*$, respectively. Note, by including the puncture balls ($r = 0$) to be annuli, we are abusing notation. The variational approach to Geometric Function Theory [3, 21, 35] adds its own demand to study such problems. Indeed, several papers are devoted to understand the expected radial symmetric properties [4, 7, 14, 16, 20, 22, 24, 25, 26, 28, 30, 31, 32, 33, 36, 37, 39]. Many times experimentally known answers to practical problems has led us to deeper insights of such mathematically challenging problems. The present paper has its origin in the radial cavitating minimizers obtained by J. Ball [7]. We seek to minimize the $p$-harmonic energy of $h$,

(1.1) $$\mathbb{E}_p[h] = \int_{\mathbb{A}} |Dh(x)|^p \, \mathrm{d}x .$$

The infimum is subjected to orientation preserving Sobolev homeomorphisms $h \colon \mathbb{A} \xrightarrow{\text{onto}} \mathbb{A}^*$ in $\mathscr{W}^{1,p}(\mathbb{A}, \mathbb{R}^n)$ which are furthermore assumed to preserve the order of the

---

2010 *Mathematics Subject Classification*. Primary 35J60; Secondary 30C70.
A. Koski was supported by the ERC Starting Grant number 307023. J. Onninen was supported by the NSF grant DMS-1700274.

boundary components,

$$
(1.2) \qquad \begin{cases} |f(x)| \to r_* & \text{when } |x| \to r \\ |f(x)| \to R_* & \text{when } |x| \to R \,. \end{cases}
$$

Such a class of Sobolev homeomorphisms is denoted by $\mathcal{H}^{1,p}(\mathbb{A}, \mathbb{A}^*)$. Here and in what follows we use the Hilbert-Schmidt norm of the differential matrix, $|Dh|^2 = \langle Dh, Dh \rangle = \operatorname{Tr}[D^*h \cdot Dh]$. Our primary question is whether the mappings of least $p$-harmonic energy are (modulo rotation) radially symmetric mappings. The question is motivated by the $2D$-models. More accurately, let $\mathbb{A}, \mathbb{A}^* \subset \mathbb{R}^2$ does the equality

$$
(1.3) \qquad \inf_{\mathcal{H}^{1,p}(\mathbb{A}, \mathbb{A}^*)} \mathbb{E}_p[h] = \inf_{\mathcal{R}^{1,p}(\mathbb{A}, \mathbb{A}^*)} \mathbb{E}_p[h]
$$

hold? In what follows, we denote the subclass of radial homeomorphisms by

$$
\mathcal{R}^{1,p}(\mathbb{A}, \mathbb{A}^*) = \left\{ h \in \mathcal{H}^{1,p}(\mathbb{A}, \mathbb{A}^*) \colon h(x) = H(|x|) \frac{x}{|x|} \right\}.
$$

When $\mathbb{A} = \mathbb{A}^* = A(0, 1) \subset \mathbb{R}^n$ we have the identity mapping $\operatorname{Id}(x) = x$ from $\mathbb{A}$ onto $\mathbb{A}^*$. Showing that the identity mapping is a minimizer for the $p$-harmonic energy provided $p \geqslant n$ is easy. Indeed, if $h \in \mathcal{H}^{1,p}(\mathbb{A}, \mathbb{A}^*)$, we can estimate the energy of $h$ using Hölder's inequality and Hadamard's inequality $|Dh|^n \geqslant n^{\frac{n}{2}} \det Dh$ as follows

$$
(1.4) \qquad \begin{aligned} \mathbb{E}_p[h] &= \int_{\mathbb{A}} |Dh(x)|^p \, \mathrm{d}x \geqslant |\mathbb{A}|^{-\frac{p-n}{n}} \left( \int_{\mathbb{A}} |Dh(x)|^n \, \mathrm{d}x \right)^{\frac{p}{n}} \\ &\geqslant |\mathbb{A}|^{-\frac{p-n}{n}} \left( n^{\frac{n}{2}} \int_{\mathbb{A}} \det Dh(x) \, \mathrm{d}x \right)^{\frac{p}{n}} = n^{\frac{p}{2}} |\mathbb{A}| = \mathbb{E}_p[\operatorname{Id}] \,. \end{aligned}
$$

When $n = 2$, it has been proven that the radial $p$-harmonic minimizers are always absolute minimizers if $p \geqslant n$, see [4, 15, 28].

**Proposition 1.1.** *Suppose that $\mathbb{A}, \mathbb{A}^* \subset \mathbb{R}^2$. If $p \geqslant 2$, then the equality* (1.3) *always holds.*

The proof of Propostion 1.1 is based on the method of *Free-Lagrangians*. A free Lagrangian (which narrows a bit the notion of null Lagrangian [5]) is a nonlinear differential $n$-form $\mathbf{L}(x, h, Dh) \, \mathrm{d}x$ defined on Sobolev homeomorphisms $h \colon \mathbb{A} \xrightarrow{\text{onto}} \mathbb{A}^*$ whose integral depends only on the homotopy class of $h$, see [28]. The volume form is not only a trivial example of free Lagrangian but also a key player in the proof. Recall that our motivation comes from the radial cavity deformations which have the infimum energy in $\mathcal{R}^{1,p}(\mathbb{A}, \mathbb{A}^*)$. In the model configuration one takes $\mathbb{A} = A(0, R)$ and $\mathbb{A}^* = A(1, R_*)$. In such model case the class of admissible homeomorphisms, $\mathcal{R}^{1,p}(\mathbb{A}, \mathbb{A}^*)$, is nonempty only when $p < n$. There are many reasons why the case $p < n$ differs from the known results. On the one hand, it is exactly the unavailability of the volume form as a free Lagrangian that causes the main mathematical challenge; new tools are needed. On the other hand, already in the planar case there is a new and unexpected phenomena concerning the failure of radial symmetry of these extremal mappings, Theorem 1.2.

1.1. **Free boundary value problem.** It is necessary to emphasize that we are not dealing with the classical $p$-harmonic boundary value problem:

- the boundary data is not given
- admissible deformations are homeomorphisms, $h(\mathbb{A}) = \mathbb{A}^*$.

The first condition, known as *traction free problems* [5, 8, 9] in the theory of nonlinear hyperelasticity, allows for "tangential slipping" along the boundary. Our first result states that in opposite of (1.4), if $p < 2$ the identity mapping $\mathrm{Id}(x) = x$ from the punctured disk onto itself need not be an absolute minimizer for the $p$-harmonic energy.

**Theorem 1.2.** *Let $\mathbb{A} = A(0,1) = \mathbb{A}^*$ be the planar punctured unit disk. Then*

$$(1.5) \quad \inf_{\mathcal{H}^{1,1}(\mathbb{A},\mathbb{A}^*)} \int_{\mathbb{A}} |Dh(x)|\,\mathrm{d}x \leqslant 2 < \sqrt{2}\pi = \mathbb{E}_1[\mathrm{Id}] = \inf_{\mathscr{R}^{1,1}(\mathbb{A},\mathbb{A}^*)} \int_{\mathbb{A}} |Dh(x)|\,\mathrm{d}x\,.$$

*Furthermore, there exists $p_1 \in (1,2)$ such that for $1 \leqslant p < p_1$ we have*

$$(1.6) \quad \inf_{\mathcal{H}^{1,p}(\mathbb{A},\mathbb{A}^*)} \mathbb{E}_p[h] < \inf_{\mathscr{R}^{1,p}(\mathbb{A},\mathbb{A}^*)} \mathbb{E}_p[h]\,.$$

It is somewhat surprising too that the infimum energy cannot be achieved within *spherically symmetric* mappings. A homeomorphism $h\colon \mathbb{A} \xrightarrow{\text{onto}} \mathbb{C}$ is spherically symmetric if it takes concentric spheres into concentric spheres; that is, $|h(x)|$ depends only on $|x|$.

1.2. **Fixed boundary value problem.** Our main purpose in this paper is to study the equality (1.3) under fixed boundary values. In view of previous results the behavior of our admissible maps on the inner boundary plays less of a role in the radial symmetry question. Furthermore keeping the model configuration $h\colon A(0,R) \xrightarrow{\text{onto}} A(1,R_*)$ in mind, we are led to impose a partial boundary data on our admissible maps. Indeed, we fix the homeomorphisms on the outer boundary of $\mathbb{A}$, allowing them to be traction free on the inner boundary. For simplifying the notation, we write $\partial_\circ \mathbb{A} = \{x \in \mathbb{R}^n \colon |x| = R\}$ and

$$\mathcal{H}^{1,p}_{\mathrm{Id}}(\mathbb{A}, \mathbb{A}^*) = \{h \in \mathscr{H}^{1,p}(\mathbb{A}, \mathbb{A}^*) \colon h \text{ is continuous up to } \partial_\circ \mathbb{A} \text{ and } h(x) = \frac{R_*}{R}x\}\,.$$

We prove that keeping the homeomorphisms fixed (non traction free) in the minimizing sequence on the outer boundary leads the hunted radial symmetry property.

**Theorem 1.3.** *[Fixed boundary values] Let $\mathbb{A}$ and $\mathbb{A}^*$ be planar annuli and $1 \leqslant p < 2$. Then*

$$(1.7) \quad \inf_{\mathcal{H}^{1,p}_{\mathrm{Id}}(\mathbb{A},\mathbb{A}^*)} \int_{\mathbb{A}} |Dh(x)|^p\,\mathrm{d}x = \inf_{\mathscr{R}^{1,p}(\mathbb{A},\mathbb{A}^*)} \int_{\mathbb{A}} |Dh(x)|^p\,\mathrm{d}x\,.$$

It turns out that the weak limit of the energy minimizing sequence of homeomorphisms need not be homeomorphisms. In fact a part of the domain near its boundary may collapse into the boundary of the target domain. In mathematical models of nonlinear elasticity this is interpreted as *interpenetration of matter* [6, 40]. We call such occurrence the *Nitsche phenomenon*, after Nitsche's conjecture [34] about existence of harmonic homeomorphisms between annuli. Now, a theorem [22]. A similar collapsing phenomenon occurs in free boundary problems for minimal graphs, called *edge-creeping* [12, 19, 41]. To discuss the Nitsche phenomenon in Theorem 1.3 we need first to enlarge the class of radially symmetric homeomorphisms. This leads us to define

$$\overline{\mathcal{R}}^{1,p}(\mathbb{A}, \mathbb{A}^*) = \{h \in \mathscr{W}^{1,p}(\mathbb{A}, \mathbb{R}^n) \colon h(x) = H(|x|) \frac{x}{|x|} \text{ and } H \in \mathscr{R}\}$$

where

$$\mathscr{R} = \{H \colon [r, R] \to [r_*, R_*] \colon H(r) = r_*, H(R) = R_* \text{ and } \dot{H} \geqslant 0\}.$$

The following result characterizes the existence of energy-minimal radially symmetric homeomorphisms. Our existence result relies only on the conformal moduli of the annuli. Of course, the $p$-energy minimization problem ($p \neq 2$) is neither invariant with respect to a conformal change of variable in the domain, nor in the target.

**Theorem 1.4** (The radial case). *Let $\mathbb{A}$ and $\mathbb{A}^*$ be planar annuli and $1 \leqslant p < 2$. There exists an increasing function $m_p \colon (1, \infty] \to (1, \infty]$ such that*

(1) *if $1 < p < 2$, then $\mathbb{E}_p$ has a unique minimizer in the class $\overline{\mathcal{R}}^{1,p}(\mathbb{A}, \mathbb{A}^*)$. Moreover, this minimizer lies in the class $\mathcal{R}^{1,p}(\mathbb{A}, \mathbb{A}^*)$ (i.e. is a homeomorphism) whenever $m_p(R/r) \leqslant R_*/r_*$.*
(2) *if $p = 1$, then the existence of a minimizer for $\mathbb{E}_1$ depends on which of the following inequalities holds*

$$\begin{cases} R_*/r_* < m_1(R/r), \text{ a unique minimizer exists in } \overline{\mathcal{R}}^{1,1}(\mathbb{A}, \mathbb{A}^*) \setminus \mathcal{R}^{1,1}(\mathbb{A}, \mathbb{A}^*). \\ m_1(R/r) \leqslant R_*/r_* \leqslant m_1^{-1}(R/r), \text{ a unique minimizer exists in } \mathcal{R}^{1,1}(\mathbb{A}, \mathbb{A}^*). \\ m_1^{-1}(R/r) < R_*/r_*, \text{ no minimizer exists in } \overline{\mathcal{R}}^{1,1}(\mathbb{A}, \mathbb{A}^*). \end{cases}$$

Specifically, the function $m_p \colon (1, \infty] \to (1, \infty]$ is given by the formula (2.10). To formulate Theorem 1.3 in its full generality, we need to adopt the limits of Sobolev homeomorphisms as legitimate deformations [29]. We denote the class of strong $\mathscr{W}^{1,p}$-limits of sequences in $\mathcal{H}^{1,p}_{\text{Id}}(\mathbb{A}, \mathbb{A}^*)$ by $\overline{\mathcal{H}}^{1,p}_{\text{Id}}(\mathbb{A}, \mathbb{A}^*)$

**Theorem 1.5.** *[Fixed boundary values] Let $\mathbb{A}$ and $\mathbb{A}^*$ be planar annuli and $1 \leq p < 2$. Then except for the case when $p = 1$ and $m_1^{-1}(R/r) < R_*/r_*$, there exists a radially symmetric mapping $h_\circ \in \overline{\mathcal{R}}^{1,p}(\mathbb{A}, \mathbb{A}^*)$ such that*

$$\min_{\overline{\mathcal{H}}^{1,p}_{\text{Id}}(\mathbb{A}, \mathbb{A}^*)} \mathbb{E}_p[h] = \min_{\overline{\mathcal{R}}^{1,p}(\mathbb{A}, \mathbb{A}^*)} \mathbb{E}_p[h] = \mathbb{E}_p[h_\circ]. \tag{1.8}$$

*The map $h_0$ is the unique minimizer in the class $\overline{\mathcal{H}}^{1,p}_{\text{Id}}(\mathbb{A}, \mathbb{A}^*)$. Furthermore, the minimizer $h_\circ$ is a homeomorphism if and only if $m_p(R/r) \leqslant R_*/r_*$.*

If $m_p(R/r) \leqslant R_*/r_*$, then the minimizer turns out to be diffeomorphism. Any diffeomorphic minimizer satisfies the *p-harmonic system*

$$\text{div}(|Dh|^{p-2} Dh) = 0, \tag{1.9}$$

since one can perform first variations $h + \varepsilon \varphi$ while preserving the diffeomorphism property. In general, however, the $p$-harmonic system (1.9) is unavailable; one cannot perform first variations within the class of Sobolev homeomorphisms (or even within the class of Sobolev mappings with nonnegative Jacobian determinants). The unavailability of the Euler-Lagrange equation is a major source of difficulty here. Such a difficulty is well recognized in the theory of nonlinear elasticity [10, 11, 38]. This is why we can only rely on the inner variation of the independent variable $x_\varepsilon = x + \varepsilon\tau(x)$, which leads to

$$(1.10) \qquad \frac{\partial}{\partial \bar{z}}\left[|Dh|^{p-2}h_z\overline{h_{\bar{z}}}\right] = \frac{2-p}{4p}\frac{\partial}{\partial z}|Dh|^p$$

Here $\frac{\partial h}{\partial z} = h_z$ and $\frac{\partial h}{\partial \bar{z}} = h_{\bar{z}}$ are complex partial derivatives of $h$. Our argument, however, does not make a direct use of the inner variational equation (1.10). The first step in the proof is to reduce the problem of minimizing the $p$-harmonic energy to a specific weighted 1-energy. The weight depends on the independent variable.

**Remark 1.6.** By virtue of the density of diffeomorphisms in $\mathcal{H}^{1,p}(\mathbb{A},\mathbb{A}^*)$, see [18, 23], we can equivalently replace the admissible homeomorphisms in (1.3) by sense-preserving diffeomorphims. Indeed, for $p \geqslant 1$, we have

$$(1.11) \qquad \inf_{\mathcal{H}^{1,p}(\mathbb{A},\mathbb{A}^*)} \mathbb{E}_p[h] = \inf_{\mathrm{Diff}(\mathbb{A},\mathbb{A}^*)} \mathbb{E}_p[h]\,.$$

Here the notation $\mathrm{Diff}(\mathbb{A},\mathbb{A}^*)$ stands for the class of orientation preserving diffeomorphisms from $\mathbb{A}$ onto $\mathbb{A}^*$ which also preserve the order of the boundary components. A similar result with fixed outer boundary $\partial_\circ \mathbb{A}$ values reads as

$$(1.12) \qquad \inf_{\mathcal{H}^{1,p}_{\mathrm{Id}}(\mathbb{A},\mathbb{A}^*)} \mathbb{E}_p[h] = \inf_{\mathrm{Diff}_{\mathrm{Id}}(\mathbb{A},\mathbb{A}^*)} \mathbb{E}_p[h]\,.$$

Here,

$$\mathrm{Diff}_{\mathrm{Id}}(\mathbb{A},\mathbb{A}^*) = \{h \in \mathrm{Diff}(\mathbb{A},\mathbb{A}^*) \colon h \text{ is continuous up to } \partial_\circ\mathbb{A} \text{ and } h(x) = \frac{R_*}{R}x\}\,.$$

In the traction free setting there is an extra symmetry for the existence of homeomorphic minimizer if $p = 1$. Indeed, the minimizer is homeomorphism provided the target annulus is neither conformally too thick nor too thin compare to the domain annulus $\mathbb{A}$. If $p > 1$ the only thinnest of $\mathbb{A}^*$ is an obstacle for the existence of homeomorphic minimizers. The next remark explains the additional symmetry in the case $p = 1$.

**Remark 1.7.** The change of variables formula, valid for for diffeomorphisms $h \colon \mathbb{A} \xrightarrow{\mathrm{onto}} \mathbb{A}^*$, relates the $\mathbb{E}_1$-energy of $h$ to the energy of its inverse as follows:

$$(1.13) \qquad \int_{\mathbb{A}} |Dh(x)|\,\mathrm{d}x = \int_{\mathbb{A}^*} |Dh^{-1}(y)|\,\mathrm{d}y\,.$$

In general the inverse of a Sobolev homeomorphism need not be a Sobolev mapping. Assuming that a Sobolev homeomorphism $h \colon \Omega \xrightarrow{\mathrm{onto}} \Omega'$ has so called *finite distortion* is the least regularity to imply that the inverse $h^{-1}$ is a Sobolev mapping [17]. Our admissible homeomorphisms are arbitrary. However, we have

$$(1.14) \qquad \inf_{h \in \mathcal{H}^{1,1}(\mathbb{A},\mathbb{A}^*)} \int_{\mathbb{A}} |Dh(x)|\,\mathrm{d}x = \inf_{h^{-1} \in \mathcal{H}^{1,1}(\mathbb{A}^*,\mathbb{A})} \int_{\mathbb{A}^*} |Dh^{-1}(x)|\,\mathrm{d}x$$

Indeed, (1.14) follows from (1.13) via (1.11). This result tells us that for $p = 1$ the minimization problem from the annulus $\mathbb{A}$ onto $\mathbb{A}^*$ is equivalent with the corresponding inverse problem from $\mathbb{A}^*$ onto $\mathbb{A}$.

We end the introduction commenting about the challenges to generalize Theorem 1.3 to the case of $n = 3$.

1.3. **3D case.** In the traction free setting there are no surprises in $\mathbb{R}^3$ when $p \geqslant 3$. It is known that there is a full analogue version of Proposition 1.1 when $n = 3$ but not when $n \geqslant 4$ [15, 28]. It was discovered in [27] that there is a stronger version of Theorem 1.2 when $n \geqslant 3$. Namely, in higher dimensions the infimum of the $n - 1$-harmonic energy among all admissible homeomorphisms is even zero.

**Proposition 1.8.** *Let $\mathbb{A} = A(0, 1) = \mathbb{A}^*$ be the punctured unit ball in $\mathbb{R}^3$. Then*

$$(1.15) \qquad \inf_{\mathcal{H}^{1,2}(\mathbb{A},\mathbb{A}^*)} \int_\mathbb{A} |Dh(x)|^2\, dx = 0 < 4\pi = \inf_{\mathscr{R}^{1,2}(\mathbb{A},\mathbb{A}^*)} \int_\mathbb{A} |Dh(x)|^2\, dx.$$

Under fixed boundary values, it is not clear for which $p$-harmonic energy one should expect to have a 3-dimensional version of Theorem 1.3. One might suspect that, by analogy, after fixing homeomorphisms on the outer boundary the deformations of smallest $p$-harmonic energy must always be radially symmetric for all $p \geqslant 1$. However, our proof of Theorem 1.3 is difficult to generalize to higher dimensions due to the lack of an analogous result to Lemma 3.1. In particular, the lemma fails due to the existence of surfaces $\phi : S^2 \to \mathbb{R}^3$ with arbitrarily small area which are also arbitrarily close to the identity map in an average (for example, $L^1$) sense, a situation which is not possible in two dimensions.

## 2. The radial case

In this section we prove Theorem 1.4.

*Proof.* Without loss of generality assume, by scaling, that $R = R_* = 1$. We have claimed that a homeomorphic minimizer among radial mappings exists when either

$$p > 1 \text{ and } m_p(1/r) \leq 1/r_* \quad \text{or} \quad p = 1 \text{ and } m_1(1/r) \leq 1/r_* \leq m_1^{-1}(1/r).$$

Let us abbreviate this condition by saying that $1/r_*$ lies in the critical interval. The following lemma now describes how this radial minimizer is obtained. In fact, the equations (2.1) and (2.2) appearing in the lemma below are reformulations of the $p$-harmonic system for radial mappings, and are equivalent with the equations obtained in [1], Chapter 4.1.

**Lemma 2.1.** *For each $p \in [1, 2)$ there exists an increasing function $m_p : (1, \infty] \to (1, \infty]$ such that the following holds true. Assuming that $1/r_*$ lies in the critical interval, there exists a function $H_0 : [r, 1] \to [r_*, 1]$ which satisfies the following conditions:*

(1) *$H_0$ is strictly increasing with $H_0(r) = r_*$ and $H_0(1) = 1$.*
(2) *$H_0$ is a solution to the separable ODE*

$$(2.1) \qquad \frac{\dot{H}_0}{H_0} = \frac{\sqrt{g(s)}}{s\sqrt{1 - g(s)}},$$

*where the function $g(s) : (r, 1) \to (0, 1)$ is a solution of the separable ODE*

$$(2.2) \qquad \dot{g}(1 - (2-p)g) = \frac{2\sqrt{g}\sqrt{1-g}(1 - pg - (2-p)\sqrt{g}\sqrt{1-g})}{s}.$$

Once we have proven the existence of $H_0$, the radial minimizer will be defined by $h_0(se^{i\theta}) = H_0(s)e^{i\theta}$.

For future reference, we also note that the equation (2.2) which was used to define $g$ is equivalent to the equation

$$\text{(2.3)} \qquad \frac{d}{ds}\left(\frac{s^{2-p}\sqrt{g(s)}}{(1-g(s))^{(p-1)/2}}\right) = s^{1-p}\frac{1-pg(s)}{(1-g(s))^{p/2}}.$$

Before we proceed to the proof of Lemma 2.1, let us see why the function $h_0$ it defines gives a minimizer of $\mathbb{E}_p$ among radial mappings. Of course, we are still in the case where $R_*$ lies in the critical interval.

*Proof of minimization - case of the critical interval.*
We begin our proof of why $h_0$ is a minimizer by writing our energy functional in polar coordinates. For any $h \in \overline{\mathcal{R}}^{1,p}(\mathbb{A}, \mathbb{A}^*)$ we have that

$$\text{(2.4)} \qquad \mathbb{E}_p[h] = \int_0^{2\pi}\int_r^1 \left(|\dot{H}|^2 + \left|\frac{H}{s}\right|^2\right)^{p/2} s\,ds d\theta.$$

Given any function $h \in \overline{\mathcal{R}}^{1,p}(\mathbb{A}, \mathbb{A}^*)$, we now estimate $\mathbb{E}_p[h]$ from below using pointwise estimates on the integrand. First, the fact that $2/p > 1$ lets us utilize the inequality

$$\text{(2.5)} \qquad (a+b)^{p/2} \geq t^{1-p/2}a^{p/2} + (1-t)^{1-p/2}b^{p/2},$$

which is valid for any nonnegative real numbers $a$ and $b$ and any $t \in [0,1]$. We use this to obtain the estimate

$$s\left(|\dot{H}|^2 + \left|\frac{H}{s}\right|^2\right)^{p/2} \geq g(s)^{1-p/2}s|\dot{H}|^p + (1-g(s))^{1-p/2}s^{1-p}H^p.$$

We now apply the arithmetic-geometric mean inequality in the following form:

$$|\dot{H}|^p \geq p\frac{g(s)^{(p-1)/2}}{s^{p-1}(1-g(s))^{(p-1)/2}}|\dot{H}|H^{p-1} - (p-1)\frac{g(s)^{p/2}}{s^p(1-g(s))^{p/2}}H^p$$

When combined with the previous estimate and simplified, this gives

$$s\left(\left|\frac{H}{s}\right|^2 + |\dot{H}|^2\right)^{p/2} \geq s^{1-p}\frac{1-pg(s)}{(1-g(s))^{p/2}}H^p + \frac{s^{2-p}\sqrt{g(s)}}{(1-g(s))^{(p-1)/2}}p|\dot{H}|H^{p-1}$$

$$\geq s^{1-p}\frac{1-pg(s)}{(1-g(s))^{p/2}}H^p + \frac{s^{2-p}\sqrt{g(s)}}{(1-g(s))^{(p-1)/2}}p\dot{H}H^{p-1}$$

$$= \frac{d}{ds}\left(\frac{s^{2-p}\sqrt{g(s)}}{(1-g(s))^{(p-1)/2}}H^p\right)$$

Here we have used formulation (2.3) of the ODE the function $g$ satisfies to obtain the last equation. We now claim that the following equation holds.

$$\text{(2.6)} \qquad \int_0^{2\pi}\int_r^1 \frac{d}{ds}\left(\frac{s^{2-p}\sqrt{g(s)}}{(1-g(s))^{(p-1)/2}}H^p\right)ds d\theta = \mathbb{E}_p[h_0].$$

This is a simple consequence of the following two facts.
  (1) The expression on the left hand side of (2.6) does not depend on $H$ as we have assumed that $H(r) = r_*$ and $H(1) = 1$.

(2) The equation (2.1), used to define $H_0$, guarantees that all of the inequalities we made in the above proof actually become equalities when $h = h_0$.

In conclusion, we have that $\mathbb{E}_p[h] \geq \mathbb{E}_p[h_0]$ for all $h \in \overline{\mathcal{R}}^{1,p}(\mathbb{A}, \mathbb{A}^*)$, which shows that $h_0$ is the minimizer in this case.

*Proof of Lemma 2.1.*
To prove Lemma 2.1 we must find, for any $p \in [1, 2)$, solutions of (2.1) and (2.2) such that $H_0(r) = r_*$ and $H_0(1) = 1$. It will be easy to define $H_0$ as soon as $g$ is chosen suitably. We will take as $H_0$ any solution of (2.1) and the condition that $H_0(1) = 1$ will be easy to guarantee by rescaling $H_0$ if necessary. The fact that $H_0$ will be strictly increasing is also a simple consequence of (2.1). We may now compute that if $H_0$ solves (2.1), then

$$\int_r^1 \frac{\sqrt{g(s)}}{s\sqrt{1-g(s)}} ds = \int_r^1 \frac{\dot{H}_0}{H_0} ds = \log H_0(1) - \log H_0(r) = -\log H_0(r).$$

Hence to prove Lemma 2.1 it remains to find a solution of (2.2) such that

(2.7) $$\int_r^1 \frac{\sqrt{g(s)}}{s\sqrt{1-g(s)}} ds = -\log r_*.$$

Let us begin by investigating the general solvability of (2.2). For $t \in (0, 1)$, denote

$$\phi_p(t) = \frac{1 - (2-p)t}{2\sqrt{t}\sqrt{1-t}(1 - pt - (2-p)\sqrt{t}\sqrt{1-t})}$$

so that (2.2) is equivalent with the equation

(2.8) $$\dot{g}\phi_p(g) = 1/s$$

although with one exception: The denominator of $\phi_p$ vanishes exactly at $t = 1/2$. This shows that the constant function $g \equiv 1/2$ is always a solution of (2.2), and it is easy to verify from (2.7) that this corresponds to the case $r = r_*$ and that $h_0$ is the identity mapping in this case. Then by the Picard-Lindelöf theorem for ODE's it follows that for any other solution $g$ we will always have either $g < 1/2$ or $g > 1/2$ everywhere. We now split into two cases:

*Case 1. The target is thinner, $r_* > r$*
This case corresponds to finding a solution $g$ of (2.8) with $g < 1/2$ everywhere, since for such a solution equation (2.7) will always imply $r_* > r$. Note that $\phi_p(t)$ is positive for $t \in (0, 1/2)$, and hence has a strictly increasing integral function $\Phi_p$ on this interval. The asymptotics of the derivative $\phi_p(t)$ at $t = 0$ and $t = 1/2$ reveal that the number $\Phi_p(0)$ is finite and that $\lim_{t \to 1/2} \Phi_p(t) = +\infty$. We may assume $\Phi_p(0) = 0$, else add a constant. Thus the inverse function $\Phi_p^{-1}$ is defined and strictly increasing on $[0, \infty)$. We may now write the general solution of (2.8) in the form

(2.9) $$g(s) = \Phi_p^{-1}(\log s + C),$$

where $s \in [1, \infty)$ and $C \geq -\log r$ is any constant. We will find that achieving (2.7) is not always possible if the target is too thin, i.e. if $r_*$ is too close to 1. Hence we will define the function $m_p(x)$ by

$$(2.10) \qquad \log m_p(1/r) = \inf_g \int_r^1 \frac{\sqrt{g(s)}}{s\sqrt{1-g(s)}}\,\mathrm{d}s,$$

where the infimum is taken over the solutions $g$ of the form (2.9). It is clear that $m_p(x)$ is defined for all $x > 1$, and from (2.9) we see that the above infimum is attained when $C = -\log r$, as $\Phi_p^{-1}$ is increasing. Thus the exact formula for $m_p$ can be calculated by

$$\log m_p(1/r) = \int_r^1 \frac{\sqrt{\Phi_p^{-1}(\log s - \log r)}}{s\sqrt{1-\Phi_p^{-1}(\log s - \log r)}}\,\mathrm{d}s$$

$$= \int_1^{1/r} \frac{\sqrt{\Phi_p^{-1}(\log s)}}{s\sqrt{1-\Phi_p^{-1}(\log s)}}\,\mathrm{d}s.$$

From the latter expression it is clear that $m_p(x)$ is increasing in $x$ and we also find that $m_p(\infty) = \infty$, which is a consequence of the factor $1/s$ inside the integral and the fact that $\Phi_p^{-1}(\log s) \to 1/2$ as $s \to \infty$. We have thus seen that the condition (2.7) cannot be attained if $1/r_* < m_p(1/r)$, but to finish Case 1 we must still show that it can be attained whenever $1/r_* \in [m_p(1/r), 1/r)$. From (2.9) one can see that when $C \to \infty$, we have $g \to 1/2$ uniformly. Hence

$$\sup_g \int_r^1 \frac{\sqrt{g(s)}}{s\sqrt{1-g(s)}}\,\mathrm{d}s = -\log r,$$

where the supremum is again taken over functions as in (2.9). Hence as $C$ goes through the values in $[-\log r, \infty)$, the right hand side of (2.7) goes through all the values in $[\log(m_p(1/r)), \log(1/r))$ as wanted.

In the above calculations we have implicitly assumed that $r \neq 0$. If the domain of definition is a punctured disc, we could have that $r = 0$. However, since $m_p(\infty) = \infty$ we find that the only case when a homeomorphic radial minimizer exists is when $r_* = 0$, in which case the minimizer is the identity map. In short terms, any target annulus with a positive inner radius is too thin compared to the punctured disc. This finishes the case $r_* > r$.

*Case 2. The target is fatter, $r_* < r$*

Analogously to Case 1, this case corresponds to finding solutions with $g > 1/2$ everywhere. Note that $\phi_p(t)$ is negative for $t \in (1/2, 1)$. In this case we can find a strictly decreasing integral function $\Psi_p(t)$ with $\Psi_p(1) = 0$ and $\lim_{t \to 1/2} \Psi_p(t) = +\infty$. Thus the inverse function $\Psi_p^{-1}$ is defined and strictly decreasing on $[0, \infty)$. Again we find solutions $g$ of the form

$$(2.11) \qquad g(s) = \Psi_p^{-1}(\log s + C),$$

where $s \in [1, \infty)$ and $C \geq -\log r$ is any constant. As in the previous case, we find that $g \to 1/2$ uniformly as $C \to \infty$. To find if there is a lower limit on the radius $r_*$ we must take $C = -\log r$ in (2.11), since this gives the largest value for the integral

$$(2.12) \qquad \int_r^1 \frac{\sqrt{g(s)}}{s\sqrt{1-g(s)}} \, ds = \int_r^1 \frac{\sqrt{\Psi_p^{-1}(\log s - \log r)}}{s\sqrt{1-\Psi_p^{-1}(\log s - \log r)}} \, ds.$$

We claim that the latter integral equals $\infty$ when $p > 1$ and $m_1^{-1}(1/r)$ when $p = 1$. Let us take first $p > 1$. Note that $\Psi_p^{-1}(\log s - \log r) \to 1$ as $s \to r$. By the mean value theorem we have

$$\begin{aligned}
1 - \Psi_p^{-1}(\log s - \log r) &= (r-s)\frac{d}{ds}\left[\Psi_p^{-1}(\log s - \log r)\right]_{s=\xi} \quad \xi \in (r,s) \\
&= (r-s)\frac{1}{\xi \phi_p(\Psi_p^{-1}(\log \xi - \log r))} \\
&\leq c(s-r)\sqrt{1 - \Psi_p^{-1}(\log \xi - \log r)} \\
&\leq c(s-r)\sqrt{1 - \Psi_p^{-1}(\log s - \log r)}
\end{aligned}$$

Where we have used the asymptotics $\phi_p(t) \approx -c(1-t)^{-1/2}$ for $t \approx 1$, valid only for $p > 1$. This gives that $\sqrt{1-\Psi_p^{-1}(\log s - \log r)} \leq c(s-r)$ for $s \approx r$, showing that the integral (2.12) diverges to $+\infty$. Thus for $p > 1$ we have a homeomorphic radial minimizer whenever $r_* < r$, including the case $r_* = 0$.

For $p = 1$ we have that $\phi_1(t) \to -1$ as $t \to 1$. This will force the integral (2.12) to be finite. The fact that the upper limit equals exactly $m_1^{-1}(1/r)$ follows essentially from the symmetry given by Remark 1.7: When $p = 1$ the radial minimizer from $\mathbb{A}(r,1)$ to $\mathbb{A}(r_*,1)$ is the inverse mapping of the radial minimizer from $\mathbb{A}(r_*,1)$ to $\mathbb{A}(r,1)$. This concludes the proof of Lemma 2.1.

*Minimization - below the critical interval.*
Let us now find the minimizer in the class $\overline{\mathcal{R}}^{1,p}(\mathbb{A},\mathbb{A}^*)$ when $R_*/r_* < m_p(R/r)$, i.e. when the target annulus is too thin. We define a radial map as follows. Let $R_0 \in (r,R)$ be the exact radius for which $m_p(R/R_0) = R_*/r_*$. Let us denote by $h_1 : \mathbb{A}(R_0, R) \to \mathbb{A}(r_*, R_*)$ the minimal radial map for $\mathbb{E}_p$, which exists and is a homeomorphism thanks to the first part of the proof. Then define the radial minimizer $h_2 : \mathbb{A}(r,R) \to \mathbb{A}(r_*, R_*)$ by

$$h_2(x) = \begin{cases} h_1(x), & \text{when } x \in \mathbb{A}(R_0, R) \\ R_0 \frac{x}{|x|}, & \text{when } x \in \mathbb{A}(r, R_0) \end{cases}$$

In geometric terms, $h_2$ collapses the annulus $\mathbb{A}(r, R_0)$ to the inner boundary of $\mathbb{A}(r_*, R_*)$.

Let now $h \in \overline{\mathcal{R}}^{1,p}(\mathbb{A},\mathbb{A}^*)$. We split our energy $\mathbb{E}_p(h)$ into two parts, one over $\mathbb{A}(R_0, R)$ and one over $\mathbb{A}(r, R_0)$, and estimate as follows. In $\mathbb{A}(R_0, R)$, we make the

Figure 1. The radial minimizer fails to be injective below the critical interval, collapsing the darker shaded part onto the inner boundary.

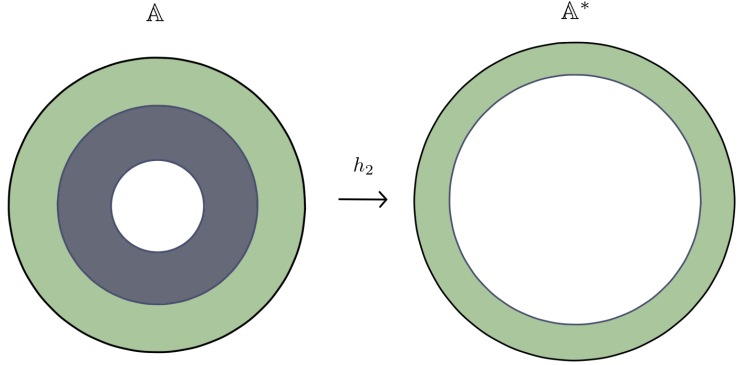

same estimates as in the case of the critical interval, ending up with

$$\int_0^{2\pi} \int_{R_0}^R \left( |\dot{H}|^2 + \left|\frac{H}{s}\right|^2 \right)^{p/2} s\,\mathrm{d}s\mathrm{d}\theta$$
$$\geq \int_0^{2\pi} \int_{R_0}^R \frac{d}{ds}\left( \frac{s^{2-p}\sqrt{g(s)}}{(1-g(s))^{(p-1)/2}} H^p \right) \mathrm{d}s\mathrm{d}\theta$$
$$= 2\pi \frac{R^{2-p}\sqrt{g(R)}}{(1-g(R))^{(p-1)/2}} R_*^p - 2\pi \frac{R_0^{2-p}\sqrt{g(R_0)}}{(1-g(R_0))^{(p-1)/2}} H(R_0)^p$$
$$= 2\pi \frac{R^{2-p}\sqrt{g(R)}}{(1-g(R))^{(p-1)/2}} R_*^p$$

Here we have used the fact that $g(R_0) = 0$. This comes from the assumption $m_p(R/R_0) = R_*/r_*$, which together with the representation (2.9) implies $g(R_0) = 0$. In $\mathbb{A}(r, R_0)$ we simply use the estimates $|\dot{H}| \geq 0$ and $H \geq 1$ to obtain

$$\int_0^{2\pi} \int_r^{R_0} \left( |\dot{H}|^2 + \left|\frac{H}{s}\right|^2 \right)^{p/2} s\,\mathrm{d}s\mathrm{d}\theta \geq 2\pi \int_r^{R_0} s^{1-p}\mathrm{d}s.$$

Adding the above two cases up gives an estimate for $\mathbb{E}_p(h)$ which is independent of $h$. One may again verify that all of the inequalities above become equalities exactly when $h = h_2$, proving our claim.

*Minimization - above the critical interval.* In this case we only need to consider what happens when $p = 1$ and $m_1^{-1}(R/r) < R_*/r_*$. However, this case will be essentially taken care of by the duality result of Remark 1.7. We have shown that for $p = 1$ and for $R/r < m_1(R_*/r_*)$, the radial minimizer from $\mathbb{A}(r_*, R_*)$ to $\mathbb{A}(r, R)$ is given by the map $h_2$ as defined above. Remark 1.7 tells us that the infimum energy $\mathbb{E}_1$ in the class $\overline{\mathcal{R}}^{1,1}(\mathbb{A}, \mathbb{A}^*)$ is exactly the $L^1$-energy of $h_2$ over $\mathbb{A}(r_*, R_*)$. However, since $h_2$ has no inverse in $\overline{\mathcal{R}}^{1,1}(\mathbb{A}, \mathbb{A}^*)$ we find that no minimal sequence for $\mathbb{E}_1$ will converge to a proper function. □

## 3. Fixed boundary value problem

In this section we describe what happens in the minimization problem subjected to homeomorphisms which fix the outer boundary, proving Theorem 1.3 and Theorem 1.5.

*Proof.* Throughout the proof, we let $h_0(re^{i\theta}) = H_0(r)e^{i\theta}$ denote the radial minimizer for $\mathbb{E}_p$ in $\overline{\mathcal{R}}^{1,p}(\mathbb{A}, \mathbb{A}^*)$. Theorem 1.4 tells us that such a mapping exists with only one exception: the case $p = 1$ and $R_*/r_*$ too large. However, this special case may be reduced back to the previous case by considering the inverse maps as in Remark 1.7. We now split into two cases:

*Case 1.* $m_p(R/r) > R_*/r_*$.

By the construction done in the previous section, there exists a radius $R_0 \in (r, R)$ so that $H_0(s) = r$ for every $s \in [r, R_0]$. Moreover, $H_0$ is defined by the equations (2.1) and (2.2) on $(R_0, R)$ and we have that $m_p(R/R_0) = R_*/r_*$. We extend the function $g : (R_0, R) \to (0, 1)$ given by Lemma 2.1 to the interval $[r, R_0]$ by setting $g(s) = 0$ for $s \in [r, R_0]$.

Let now $h \in \overline{\mathcal{H}}_{\mathrm{Id}}^{1,p}(\mathbb{A}, \mathbb{A}^*)$ be any map. In all of the following estimates equality will again hold for $h_0$ in place of $h$. We begin our estimates by utilizing Hölder's inequality in the form

$$(3.1) \qquad \int_{\mathbb{A}} |Dh|^p \mathrm{d}x \geq \frac{\left(\int_{\mathbb{A}} |Dh_0|^{p-1}|Dh|\mathrm{d}x\right)^p}{\left(\int_{\mathbb{A}} |Dh_0|^p \mathrm{d}x\right)^p}.$$

This reduces us to estimating the quantity $\int_{\mathbb{A}} |Dh_0|^{p-1}|Dh|dx$. We begin by writing $|Dh|$ in polar coordinates and using the inequality (2.5) (with 1 in place of $p$) to obtain

$$(3.2) \qquad s|Dh| = s\left(|h_s|^2 + \left|\frac{h_\theta}{s}\right|^2\right)^{1/2} \geq \sqrt{g(s)}\,s|h_s| + \sqrt{1-g(s)}\,|h_\theta|.$$

Combining, we have that

$$(3.3) \qquad \int_{\mathbb{A}} |Dh_0|^{p-1}|Dh|\,\mathrm{d}x \geq \int_0^{2\pi}\int_r^R \rho_1(s)|h_s| + \rho_2(s)|h_\theta|\,\mathrm{d}s\,\mathrm{d}\theta,$$

where

$$\rho_1(s) = s\sqrt{g(s)}\left(|\dot{H}_0|^2 + |H_0/s|^2\right)^{(p-1)/2} \quad \text{and}$$

$$\rho_2(s) = \sqrt{1-g(s)}\left(|\dot{H}_0|^2 + |H_0/s|^2\right)^{(p-1)/2}.$$

Once we are done with the rest of the estimates, the uniqueness part of Theorem 1.5 will already follow from (3.3) as there equality can only hold for $h = h_0$. This also allows us to assume that $h$ is a diffeomorphism for the rest of the proof, as to find the infimum of the energy $\mathbb{E}_p$ it is sufficient to minimize among the class $\mathrm{Diff}_{\mathrm{Id}}(\mathbb{A}, \mathbb{A}^*)$ as mentioned in Remark 1.6.

We now estimate the expression on the right hand side of (3.3) separately on the intervals $[r, R_0]$ and $(R_0, R)$:

*The interval $[r, R_0]$.*
We observe that the estimate
$$\int_0^{2\pi} |h_\theta(se^{i\theta})| \, d\theta \geq 2\pi r$$
holds for all radii $s$, as the expression on the right hand side denotes the length of the image curve of $\{se^{i\theta} : \theta \in [0, 2\pi)\}$. Thus combining with the trivial estimate $|h_s| \geq 0$, we have that
$$\int_0^{2\pi} \int_r^{R_0} \rho_1(s)|h_s| + \rho_2(s)|h_\theta| \, ds d\theta \geq 2\pi r \int_r^{R_0} \rho_2(s) ds.$$
Since on this interval one has $g(s) = 0$, we may further observe that $\rho_2(s) = s^{1-p}$ to get an explicit lower bound. However, for the proof it is simply enough to have a lower estimate independent of $h$ (and with equality for $h = h_0$).

*The interval $(R_0, R)$.*
On this interval we claim to have the curious identity $\dot\rho_1 = \rho_2$. This is simply a consequence of the equations (2.1) and (2.2) after the following computation.

$$\dot\rho_1 = \frac{d}{ds}\left[s\sqrt{g(s)}\left(|\dot H_0|^2 + |H_0/s|^2\right)^{(p-1)/2}\right]$$

$$= \frac{d}{ds}\left[s\sqrt{g(s)}\left(\left|\frac{\sqrt{g(s)}}{s\sqrt{1-g(s)}}H_0\right|^2 + |H_0/s|^2\right)^{(p-1)/2}\right]$$

$$= \frac{d}{ds}\left[\frac{s^{2-p}\sqrt{g(s)}}{(1-g(s))^{(p-1)/2}}H_0^{p-1}\right]$$

$$= H_0^{p-1}\left(s^{1-p}\frac{1-pg(s)}{(1-g(s))^{p/2}} + (p-1)\frac{\sqrt{g(s)}}{s\sqrt{1-g(s)}}\frac{s^{2-p}\sqrt{g(s)}}{(1-g(s))^{(p-1)/2}}\right)$$

$$= H_0^{p-1}\left(s^{1-p}\frac{1-g(s)}{(1-g(s))^{p/2}}\right)$$

$$= \sqrt{1-g(s)}\left(\left(\frac{g(s)}{s^2(1-g(s))} + \frac{1}{s^2}\right)H_0^2\right)^{(p-1)/2}$$

$$= \rho_2.$$

Furthermore, we have $\rho_1(R_0) = 0$ since $g(R_0) = 0$. We may now estimate that
$$\int_0^{2\pi}\int_{R_0}^R \rho_1(s)|h_s(se^{i\theta})| \, ds \, d\theta = \int_0^{2\pi}\int_{R_0}^R \int_{R_0}^s \rho_2(t)|h_s(se^{i\theta})| \, dt \, ds \, d\theta$$
$$= \int_0^{2\pi}\int_{R_0}^R \rho_2(t) \int_t^R |h_s(se^{i\theta})| \, ds \, dt \, d\theta$$
$$\geq \int_0^{2\pi}\int_{R_0}^R \rho_2(t)\left|h(Re^{i\theta}) - h(te^{i\theta})\right| \, dt \, d\theta.$$

Recalling the fact that $h(Re^{i\theta}) = R_*e^{i\theta}$ now makes the above estimate read out as
$$\int_0^{2\pi}\int_{R_0}^R \rho_1(s)|h_s| + \rho_2(s)|h_\theta| \, ds \, d\theta \geq \int_{R_0}^R \rho_2(t) \int_0^{2\pi} \left|R_*e^{i\theta} - h(te^{i\theta})\right| + |h_\theta(te^{i\theta})| \, d\theta \, dt,$$

where we have simply renamed the variable $s$ to $t$ in the second term. Notice that the quantity inside the inner integral only depends on the image curve of $\{te^{i\theta} : \theta \in [0, 2\pi)\}$ under $h$. Hence we will be done once we prove the following lemma about such curves.

**Lemma 3.1.** *For any rectifiable closed curve $\gamma : [0, 2\pi) \to \bar{\mathbb{D}}$ we have that*

$$B(\gamma) = \int_0^{2\pi} |e^{it} - \gamma(t)| + |\dot{\gamma}(t)| \, dt \geq 2\pi,$$

*with equality when $\gamma$ is any circle centered at the origin.*

*Proof of Lemma 3.1.* We begin by showing that the energy $B(\gamma)$ is convex. Take any two curves $\gamma_1$ and $\gamma_2$. Let $\gamma = (\gamma_1 + \gamma_2)/2$. We argue that $B(\gamma) \leq (B(\gamma_1) + B(\gamma_2))/2$. This follows simply from the estimates

$$|\dot{\gamma}(t)| \leq (|\dot{\gamma_1}(t)| + |\dot{\gamma_2}(t)|)/2$$

and

$$|e^{it} - \gamma(t)| \leq \left(|e^{it} - \gamma_1(t)| + |e^{it} - \gamma_2(t)|\right)/2.$$

For our second observation we note that $B(\gamma)$ is invariant under rotation of $\gamma$: For any angle $\theta \in [0, 2\pi)$, the curve $\gamma_\theta(t) = e^{i\theta}\gamma(t - \theta)$ satisfies $B(\gamma_\theta) = B(\gamma)$.

Let now $\gamma$ be any curve and $N$ be a large integer. Letting $\theta_k = 2\pi k/N$ for $k = 1, \ldots, N$, we define a new curve $\beta$ by

$$\beta(t) = \frac{1}{N} \sum_{k=1}^{N} \gamma_{\theta_k}(t).$$

From the convexity and rotational symmetry of $B$ we find that $B(\beta) \leq B(\gamma)$. We also find that the curve $\beta$ has some rotational symmetry of its own. In fact, $\beta = \beta_{\theta_k}$ for all $k$, which means that $\beta$ has $N$-fold radial symmetry. For the last part of the proof we need one more observation. We argue that if the image of $\beta$ is not the boundary of a convex region, then $\beta$ may be replaced with a curve with smaller energy. We first note that the energy $B(\beta)$ may be estimated below by

$$B(\beta) = \int_0^{2\pi} |e^{it} - \beta(t)| + \left|\dot{\beta}(t)\right| dt$$
$$\geq \int_0^{2\pi} \mathrm{dist}(e^{it}, \beta) + \left|\dot{\beta}(t)\right| dt$$

Our replacement curve $\tilde{\beta}$ will be defined as follows: It will be the boundary of the convex hull of $\beta$, parametrized (not necessarily injectively) in such a way that we have the equality

$$|e^{it} - \tilde{\beta}(t)| = \mathrm{dist}(e^{it}, \tilde{\beta}).$$

Such a parametrization always exists: For every point $\tilde{\beta}(t_0)$ on our convex curve we may draw two normal lines corresponding to the two one-sided tangents at that point. These lines will hit $\partial \mathbb{D}$ at the points $e^{it_2}$ and $e^{it_3}$. Our parametrization is then defined so that the interval $[t_2, t_3]$ gets mapped to exactly $\tilde{\beta}(t_0)$. Of course one could also have made the observation that for almost every point the curve $\tilde{\beta}$ is differentiable and at these points we will have $t_2 = t_3$.

We now see that $B(\tilde{\beta}) \leq B(\beta)$, since by taking the boundary of the convex hull we have made our curve shorter, meaning

$$\int_0^{2\pi} \left|\dot{\tilde{\beta}}(t)\right| dt \leq \int_0^{2\pi} \left|\dot{\beta}(t)\right| dt,$$

and we have also decreased the distance to the boundary, meaning

$$\operatorname{dist}(e^{it}, \tilde{\beta}) \leq \operatorname{dist}(e^{it}, \beta).$$

Now since $\beta$ is radially $N$-fold symmetric, so is the curve $\tilde{\beta}$ by construction. But as $N$ increases, any $N$-fold symmetric convex curve will be arbitrarily close to a circle centered at the origin both uniformly and in length. Since every circle has energy $2\pi$, letting $N \to \infty$ gives that $B(\gamma) \geq 2\pi$ for every curve $\gamma$. By looking at the cases of equality in the above proof, one can conclude that no other minimizers except circles centered at the origin exist. This concludes the proof of Lemma 3.1.

Combining with our earlier estimate on the interval $(R_0, R)$, we now obtain from Lemma 3.1 by scaling that

$$\int_0^{2\pi} \int_{R_0}^R \rho_1(s)|h_s| + \rho_2(s)|h_\theta| \, ds \, d\theta \geq 2\pi R_* \int_{R_0}^R \rho_2(t) \, dt.$$

If we include the estimates on the interval $[r, R_0)$ as well, we finally obtain that

$$\mathbb{E}_p[h] \geq 2\pi R_* \int_{R_0}^R \rho_2(t) \, dt + 2\pi r \int_r^{R_0} \rho_2(s) \, ds = \mathbb{E}_p[h_0].$$

This shows that $h_0$ is the minimizer.

*Case 2.* $m_p(R/r) \leq R_*/r_*$.

In this case we know that the radial minimization problem has a unique homeomorphic minimizer, which we again denote by $h_0$. We then proceed with the same estimates (3.1), (3.2) and (3.3) as in the previous case. We again find that $\rho_2 = \dot{\rho}_1$. Hence

$$\int_0^{2\pi} \int_r^R \rho_1(s)|h_s(se^{i\theta})| \, ds \, d\theta = \int_0^{2\pi} \int_r^R \int_r^s \rho_2(t)|h_s(se^{i\theta})| \, dt \, ds \, d\theta$$
$$+ \int_0^{2\pi} \int_r^R \rho_1(r)|h_s(se^{i\theta})| \, ds \, d\theta$$
$$\geq \int_0^{2\pi} \int_r^R \rho_2(t) \int_t^R |h_s(se^{i\theta})| \, ds \, dt \, d\theta$$
$$+ \int_0^{2\pi} \rho_1(r)|h(Re^{i\theta}) - h(re^{i\theta})| \, d\theta$$
$$\geq \int_0^{2\pi} \int_r^R \rho_2(t) \left|h(Re^{i\theta}) - h(te^{i\theta})\right| dt \, d\theta$$
$$+ 2\pi \rho_1(r)(R_* - r_*).$$

Compared to the previous case, there appeared an extra term due to the fact that $\rho_1(r) \neq 0$. We were able to estimate this extra term from below by the quantity $2\pi\rho_1(r)(R_* - r_*)$ which is independent of $h$ and attained for $h = h_0$. Hence the rest of the proof proceeds as in Case 1, utilizing Lemma 3.1 to show that

$$\int_0^{2\pi} \int_r^R \rho_2(t) \left(\left|h(Re^{i\theta}) - h(te^{i\theta})\right| + |h_\theta(te^{i\theta})|\right) dt\, d\theta \geq 2\pi R_* \int_r^R \rho_2(t)\, dt.$$

This shows that the unique minimizer is indeed the radial minimizer $h_0$. □

## 4. Free boundary values in the plane

In this section we prove Theorem 1.2. Throughout this section $\mathbb{D}$ denotes the open unit disc in the plane. Our first step in the proof will be to reduce to the study of a slightly simpler energy functional. Note first the elementary estimates

$$\sqrt{2}\sqrt{a^2 + b^2} \geq a + b \geq \sqrt{a^2 + b^2}.$$

These show that it will be enough to prove the following theorem.

**Theorem 4.1.** *For the energy functional*

(4.1) $$e[h] = \int_0^{2\pi} \int_0^1 s|h_s| + |h_\theta|\, ds\, d\theta$$

*we have that*

$$\inf_{h \in \mathcal{H}^{1,1}(\mathbb{A},\mathbb{A}^*)} e[h] = 2.$$

To elaborate, once we have proven Theorem 4.1 we may observe the following. For the identity mapping we have $\mathbb{E}_1[\text{Id}] = e[\text{Id}]/\sqrt{2} = \sqrt{2}\pi$. However, there exist functions $h$ such that $e[h] \approx 2$, so we have $\mathbb{E}_1[h] \leq e[h] < \mathbb{E}_1[\text{Id}]$.

*Proof of Theorem 1.2.* We construct a sequence of mappings whose energies converge to the desired limit. Let $\epsilon > 0$. We construct a mapping $h^{(\epsilon)}$ as the composition of two homeomorphisms of $\mathbb{D} \setminus \{0\}$ to itself. The mapping $h^{(\epsilon)}$ will take everyting except a small set in $\mathbb{D} \setminus \{0\}$ close to the point $(1,0)$.

*The first map.*
Our first map $H_1$ will be a mapping that keeps all of the circles $C_s = \{se^{i\theta}, 0 \leq \theta < 2\pi\}$ fixed for $0 < s < 1$. In fact, it will simply be a reparametrization of each such circle. Let $S_\epsilon(\theta_0)$ denote an interval of length $\epsilon$ around the point $\theta_0 \in \mathbb{R}$. We define our reparametrization so that it stretches the arc $A_1 = \{se^{i\theta}, \theta \in S_\epsilon(\pi)\}$ onto the complement of the arc $A_2 = \{se^{i\theta}, \theta \in S_\epsilon(0)\}$ on our circle $C_s$ and vice versa, squeezing the complement of $A_1$ into $A_2$. We let this parametrization be independent of $s$ so that the image of any radial line $R_\theta = \{se^{i\theta} : 0 < s < 1\}$ under $H_1$ is another radial line.

*The second map.*
Our second map $H_2$ will keep the radial lines $R_\theta$ fixed. How it maps each radial line will depend on the angle $t$:

(1) If $\theta \in S_\epsilon(0)$, our map will map the segment $\{se^{i\theta} : 0 < s < \epsilon\}$ linearly to the segment $\{se^{i\theta} : 0 < s < 1-\epsilon\}$ and do the same in their complementary segments on $R_\theta$. Hence it squeezes everything except a small set on $R_\theta$ close to the boundary.

(2) If $t \in (-\pi, \pi] \setminus S_{2\epsilon}(0)$, our map will behave exactly oppositely as in the previous case, stretching the segment $\{se^{i\theta} : 1 - \epsilon < s < 1\}$ linearly to $\{se^{i\theta} : \epsilon < s < 1\}$. Hence it squeezes everything except a small set on $R_\theta$ close to zero.
(3) In between the two above sets, simply interpolate the two reparametrizations linearly.

We may then define $h^{(\epsilon)} = H_2 \circ H_1$. By construction $h^{(\epsilon)}$ maps radial lines to other radial lines. The energy $e[h^{(\epsilon)}]$ will now be computed in two parts.

*The angular part.*
Here we estimate the term involving $|h_\theta^{(\epsilon)}|$ appearing in $e[h^{(\epsilon)}]$. We observe that for each $s \in (0, 1)$ the integral $\int_0^{2\pi} \left|h_\theta^{(\epsilon)}(se^{i\theta})\right| d\theta$ is equal to the length of the image curve of the circle $C_s$ under $h^{(\epsilon)}$. Since length is independent of parametrization, the mapping $H_1$ does not play a role here. We now look at the lengths of these circles under $H_2$.

FIGURE 2. Image of a circle $C_s = \{z : |z| = s\}$ under the map $H_2$.

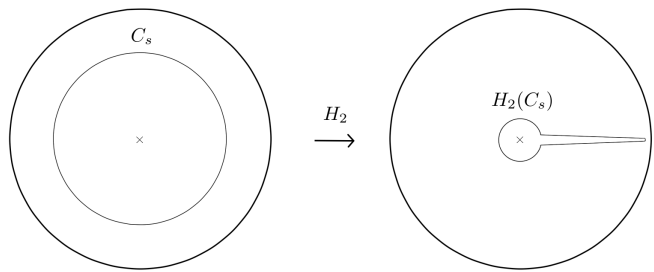

For $0 < s < 1 - \epsilon$, Figure 2 shows that the length of the image curve quickly approaches 2 as $\epsilon \to 0$. For $1 - \epsilon \leq s < 1$ the length varies but stays bounded by a constant independent of $\epsilon$. This shows that

$$(4.2) \qquad \int_0^1 \int_0^{2\pi} \left|h_\theta^{(\epsilon)}(se^{i\theta})\right| d\theta \, ds = 2 + O(\epsilon).$$

*The radial part.*
In this part we estimate the term involving $s|h_s^{(\epsilon)}|$ in $e\left[h^{(\epsilon)}\right]$. Hence for each angle $\theta$, we will estimate the value of the integral $\int_0^1 s\left|h_s^{(\epsilon)}(se^{i\theta})\right| d\theta$. Note first the trivial estimate

$$\int_0^1 s\left|h_s^{(\epsilon)}(se^{i\theta})\right| d\theta \leq \int_0^1 \left|h_s^{(\epsilon)}(se^{i\theta})\right| d\theta = 1,$$

where we have used the fact that $h^{(\epsilon)}$ maps radial lines to other radial lines to obtain the last equality. We will apply this estimate for angles in $S_\epsilon(\pi)$ - a small interval around the point $\theta = \pi$. For the other angles we argue as follows. If $\theta \notin S_\epsilon(\pi) + 2\pi\mathbb{Z}$, then the map $H_1$ first maps the radial line $R_\theta$ to a radial line $R_\alpha$ with $\alpha \in S_\epsilon(0)$. Then the map $H_2$ keeps this radial line fixed, but stretches the line in such a way that

$$\frac{d}{ds}H_2(se^{i\alpha}) = \begin{cases} \frac{1-\epsilon}{\epsilon} \text{ for } s \in (0,\epsilon) \\ \frac{\epsilon}{1-\epsilon} \text{ for } s \in (\epsilon,1) \end{cases}$$

Thus for angles $\theta \notin S_\epsilon(\pi) + 2\pi\mathbb{Z}$, we have

$$\int_0^1 s \left|h_s^{(\epsilon)}(se^{i\theta})\right| \, \mathrm{d}s = \int_0^\epsilon s \frac{1-\epsilon}{\epsilon} \, \mathrm{d}s + \int_\epsilon^1 s \frac{\epsilon}{1-\epsilon} \, \mathrm{d}s = O(\epsilon).$$

Combining this with the estimates for $\theta \in S_\epsilon(\pi)$ and with (4.2), we finally obtain that

$$e\left[h^{(\epsilon)}\right] = \int_0^{2\pi} \int_0^1 s \left|h_s^{(\epsilon)}\right| \, \mathrm{d}\theta \, \mathrm{d}s + \int_0^1 \int_0^{2\pi} \left|h_\theta^{(\epsilon)}\right| \, \mathrm{d}s \, \mathrm{d}\theta = 2 + O(\epsilon).$$

$\square$

## References

1. T. Adamowicz, *On the geometry of p-harmonic mappings*, Ph. D. thesis, (2008) Syracuse University.
2. S. S. Antman, *Nonlinear problems of elasticity. Applied Mathematical Sciences*, 107. Springer-Verlag, New York, 1995.
3. K. Astala, T. Iwaniec, and G. Martin, *Elliptic partial differential equations and quasiconformal mappings in the plane*, Princeton University Press, 2009.
4. K. Astala, T. Iwaniec, and G. Martin, *Deformations of annuli with smallest mean distortion*, Arch. Ration. Mech. Anal. **195** (2010), no. 3, 899–921.
5. J. M. Ball, *Convexity conditions and existence theorems in nonlinear elasticity*, Arch. Rational Mech. Anal. **63** (1976/77), no. 4, 337–403.
6. Ball, J.M. *Global invertibility of Sobolev functions and the interpenetration of matter.* Proc. Roy. Soc. Edinburgh Sect. A 88 (1981), no. 3-4, 315–328.
7. J. M. Ball, *Discontinuous equilibrium solutions and cavitation in nonlinear elasticity*, Philos. Trans. R. Soc. Lond. A **306** (1982) 557–611.
8. J. M. Ball, *Constitutive inequalities and existence theorems in nonlinear elastostatics*, Nonlinear analysis and mechanics: Heriot-Watt Symposium (Edinburgh, 1976), Vol. I, pp. 187–241. Res. Notes in Math., No. 17, Pitman, London, (1977).
9. J. M. Ball, *Existence of solutions in finite elasticity*, Proceedings of the IUTAM Symposium on Finite Elasticity. Martinus Nijhoff, 1981.
10. J. M. Ball, *Minimizers and the Euler-Lagrange equations*, Trends and applications of pure mathematics to mechanics (Palaiseau, 1983), 1–4, Lecture Notes in Phys., 195, Springer, Berlin, 1984.
11. J. M. Ball, *Some open problems in elasticity*, Geometry, mechanics, and dynamics, 3–59, Springer, New York, 2002.
12. Y. W. Chen, *Discontinuity and representations of minimal surface solutions*, Proceedings of the conference on differential equations (dedicated to A. Weinstein), pp. 115–138. University of Maryland, College Park, MD (1956).
13. P. G. Ciarlet, *Mathematical elasticity Vol. I. Three-dimensional elasticity*, Studies in Mathematics and its Applications, 20. North-Holland Publishing Co., Amsterdam, 1988.
14. J.-M. Coron, R.D. Gulliver, *Minimizing p-harmonic maps into spheres*, J. Reine Angew. Math. **401** (1989) 82–100.
15. D. Cuneo, *Mappings between annuli of mimimal p-Harmonic energy*, Ph. D. thesis, (2017) Syracuse University.
16. R. Hardt, F.H. Lin, C.Y. Wang, *The p-energy minimality of $x/|x|$*, Comm. Anal. Geom. 6 (1998) 141–152
17. S. Hencl and P. Koskela, *Regularity of the inverse of a planar Sobolev homeomorphism.* Arch. Ration. Mech. Anal. 1**80** (2006), no. 1, 75–95.
18. S. Hencl and A. Pratelli *Diffeomorphic Approximation of $W^{1,1}$ Planar Sobolev Homeomorphisms*, J. Eur. Math. Soc. to appear.

19. S. Hildebrandt and J. C. C. Nitsche, *A uniqueness theorem for surfaces of least area with partially free boundaries on obstacles*, Arch. Rational Mech. Anal. **79** (1982), no. 3, 189–218.
20. M.-C. Hong, *On the minimality of the p-harmonic map $\frac{x}{|x|} \colon B^n \to S^{n-1}$*, Calc. Var. Partial Differential Equations 13 (2001) 459–468.
21. T. Iwaniec and G. Martin, *Geometric Function Theory and Non-linear Analysis*, Oxford Mathematical Monographs, Oxford University Press, 2001.
22. T. Iwaniec, L. V. Kovalev, and J. Onninen, *The Nitsche conjecture*, J. Amer. Math. Soc. **24** (2011), no. 2, 345–373.
23. T. Iwaniec, L. V. Kovalev, and J. Onninen, *Diffeomorphic approximation of Sobolev homeomorphisms*, Arch. Rat. Mech. Anal. **201** (2011), no. 3, 1047–1067.
24. T. Iwaniec, L. V. Kovalev, and J. Onninen, *Doubly connected minimal surfaces and extremal harmonic mappings*, J. Geom. Anal. **22** (2012), no. 3, 726–762.
25. T. Iwaniec and J. Onninen, *Hyperelastic deformations of smallest total energy*, Arch. Ration. Mech. Anal. **194** (2009), no. 3, 927–986.
26. T. Iwaniec and J. Onninen, *Neohookean deformations of annuli, existence, uniqueness and radial symmetry*, Math. Ann. **348** (2010), no. 1, 35–55.
27. T. Iwaniec and J. Onninen, *An invitation to n-harmonic hyperelasticity*, Pure Appl. Math., Q., **7**, (2011), Special Issue: In honor of Frederick W. Gehring, Part 2, 319 –343.
28. T. Iwaniec and J. Onninen, *n-Harmonic mappings between annuli* Mem. Amer. Math. Soc. **218** (2012).
29. T. Iwaniec and J. Onninen, *Limits of Sobolev homeomorphisms* J. Eur. Math. Soc. (JEMS) **19** (2017), no. 2, 473–505.
30. M. Jordens and G. J. Martin, *Deformations with smallest weighted $L^p$ average distortion and Nitsche type phenomena*, J. Lond. Math. Soc. (2) **85** (2012), no. 2, 282–300.
31. W. Jäger, H. Kaul, *Rotationally symmetric harmonic maps from a ball into a sphere and the regularity problem for weak solutions of elliptic systems* J. Reine Angew. Math. 343 (1983) 146–161.
32. F. Meynard, *Existence and nonexistence results on the radially symmetric cavitation problem*, Quart. Appl. Math. **50** (1992) 201–226.
33. S. Müller, S. J. Spector, *An existence theory for nonlinear elasticity that allows for cavitation*, Arch. Rational Mech. Anal. **131** (1995) 1–66.
34. J.C.C. Nitsche, *On the modulus of doubly connected regions under harmonic mappings*, Amer. Math. Monthly, **69**, (1962), 781–782.
35. Yu. G. Reshetnyak, *Space mappings with bounded distortion*, American Mathematical Society, Providence, RI, 1989.
36. J. Sivaloganathan, *Uniqueness of regular and singular equilibria for spherically symmetric problems of nonlinear elasticity*, Arch. Rational Mech. Anal. 96 (1986) 97–136.
37. J. Sivaloganathan and S. J. Spector, *Necessary conditions for a minimum at a radial cavitating singularity in nonlinear elasticity*, Ann. Inst. H. Poincaré Anal. Non Linéaire **25** (2008), no. 1, 201–213.
38. J. Sivaloganathana and S. J. Spector, *On irregular weak solutions of the energy-momentum equations*, Proc. R. Soc. Edinb. A **141** (2011), 193–204.
39. C.A. Stuart, *Radially symmetric cavitation for hyperelastic materials*, Anal. Non Liné aire 2 (1985) 33?66.
40. V. Šverák, *Regularity properties of deformations with finite energy*, Arch. Rational Mech. Anal. **100** (1988), no. 2, 105–127.
41. G. Turowski, *Behaviour of doubly connected minimal surfaces at the edges of the support surface*, Arch. Math. (Basel) **77** (2001), no. 3, 278–288.

Department of Mathematics and Statistics, P.O.Box 35 (MaD) FI-40014 University of Jyväskylä, Finland
*E-mail address*: `aleksis.t.koski@jyu.fi`

Department of Mathematics, Syracuse University, Syracuse, NY 13244, USA and Department of Mathematics and Statistics, P.O.Box 35 (MaD) FI-40014 University of Jyväskylä, Finland
*E-mail address*: `jkonnine@syr.edu`